\pdfoutput=1
\RequirePackage{ifpdf}
\ifpdf 
\documentclass[pdftex]{sigma}
\else
\documentclass{sigma}
\fi

\renewcommand{\d}{\mathrm{d}}
\newcommand{\cL}{{\mathcal L}}
\newcommand{\cI}{{\mathcal I}}
\newcommand{\bbR}{{\mathbb R}}

\newcommand{\bbP}{{\mathbb P}}

\newcommand{\Gr}{\operatorname{Gr}}
\newcommand{\Leg}{\operatorname{Leg}}
\newcommand{\Lag}{\operatorname{Lag}}
\newcommand{\Ss}{\mathsf{S}}
\newcommand{\w}{\wedge}

\numberwithin{equation}{section}

\newtheorem{Theorem}{Theorem}[section]
\newtheorem*{Theorem*}{Theorem}
\newtheorem{Corollary}[Theorem]{Corollary}

\newtheorem{Proposition}[Theorem]{Proposition}
 { \theoremstyle{definition}
\newtheorem{Definition}[Theorem]{Definition}

\newtheorem{Example}[Theorem]{Example}
\newtheorem{Remark}[Theorem]{Remark} }

\begin{document}

\allowdisplaybreaks

\renewcommand{\thefootnote}{}

\newcommand{\arXivNumber}{1512.07100}

\renewcommand{\PaperNumber}{060}

\FirstPageHeading

\ShortArticleName{On the Convex Pfaff--Darboux Theorem of Ekeland and Nirenberg}

\ArticleName{On the Convex Pfaff--Darboux Theorem\\ of Ekeland and Nirenberg\footnote{This paper is a~contribution to the Special Issue on Differential Geometry Inspired by Mathematical Physics in honor of Jean-Pierre Bourguignon for his 75th birthday. The~full collection is available at \href{https://www.emis.de/journals/SIGMA/Bourguignon.html}{https://www.emis.de/journals/SIGMA/Bourguignon.html}}}

\Author{Robert L.~BRYANT}

\AuthorNameForHeading{R.L.~Bryant}

\Address{Department of Mathematics, Duke University,\\ PO Box 90320, Durham, NC 27708-0320, USA}
\Email{\href{mailto:bryant@math.duke.edu}{bryant@math.duke.edu}}
\URLaddress{\url{https://fds.duke.edu/db/aas/math/faculty/bryant}}

\ArticleDates{Received July 20, 2023, in final form August 20, 2023; Published online August 23, 2023}

\Abstract{The classical Pfaff--Darboux theorem, which provides local `normal forms' for $1$-forms on manifolds, has applications in the theory of certain economic models [Chiappo\-ri~P.-A., Ekeland I., \textit{Found. Trends Microecon.} \textbf{5} (2009), 1--151]. However, the normal~forms needed in these models often come with an additional requirement of some type of~\emph{con\-vexity}, which is not provided by the classical proofs of the Pfaff--Darboux theorem. (The appropriate notion of `convexity' is a feature of the economic model. In the simplest case, when the economic model is formulated in a domain in~$\mathbb{R}^n$, convexity has its usual meaning.) In [\textit{Methods Appl. Anal.} \textbf{9} (2002), 329--344], Ekeland and Nirenberg were able to characterize necessary and sufficient conditions for a given $1$-form~$\omega$ to admit a convex local normal form (and to show that some earlier attempts [Chiappori~P.-A., Ekeland I., \textit{Ann. Scuola Norm.\ Sup.\ Pisa Cl.\ Sci.~4} \textbf{25} (1997), 287--297] and [Zakalyukin V.M., \textit{C.~R.~Acad.\ Sci.\ Paris S\'{e}r.~{I} Math.} \textbf{327} (1998), 633--638] at this characterization had been unsuccessful). In this article, after providing some necessary background, I prove a strengthened and generalized convex Pfaff--Darboux theorem, one that covers the case of a Legendrian foliation in which the notion of convexity is defined in terms of a torsion-free affine connection on the underlying manifold. (The main result of Ekeland and Nirenberg concerns the case in which the affine connection is flat.)}

\Keywords{Pfaff--Darboux theorem; convexity; utility theory}

\Classification{58A15; 91B16}

\renewcommand{\thefootnote}{\arabic{footnote}}
\setcounter{footnote}{0}

\section{Introduction}

The Pfaff--Darboux theorem provides a local `normal form'
for $1$-forms on manifolds, assuming that certain constant rank
conditions are met. A common version\footnote{There are many variants. See~\cite[Chapter II, Section~3]{MR1083148}.}
of this classical theorem
is the following: Let $\omega$ be a smooth $1$-form on an $n$-manifold~$M$
and suppose that there is an integer~$k>0$ such that
\[
\text{$\omega\w(\d\omega)^k$ vanishes identically on $M$}
\]
while
\[
\text{$\omega\w(\d\omega)^{k-1}$ is nowhere vanishing on $M$.}
\]
Then each $m\in M$ has an open neighborhood~$U\subset M$
on which there exist (smooth) functions~$y^1,\ldots,y^k$, $p_2,\ldots,p_k$
and a nonvanishing function~$a$ such that\footnote{Throughout this article, I adopt the convention that,
when $L\subset M$ is a submanifold and $\psi$ is a differential
form on~$M$, then $L^*\psi$ denotes the \emph{pullback} of $\psi$ to~$L$.}
\begin{equation}\label{eq: P-Dnormalform1}
U^*\omega = a (\d y^1 + p_2\,\d y^2 + \cdots + p_k\,\d y^k).
\end{equation}
Since
\[
U^*\bigl(\omega\w(\d\omega)^{k-1}\bigr)
=(-1)^{k(k-1)/2}(k{-}1)!\,a^k\,\d y^1\w\cdots\w\d y^k\w\d p_2\w\cdots\w\d p_k,
\]
the functions $y^1,\ldots,y^k$, $p_2,\ldots,p_k$
in this representation must be independent on~$U$.

The normal form~\eqref{eq: P-Dnormalform1}
is often written more symmetrically as
\begin{equation}\label{eq: P-Dnormalform2}
U^*\omega = a_1\,\d u^1 + a_2\,\d u^2 + \cdots + a_k\,\d u^k,
\end{equation}
where the $a_i$ do not simultaneously vanish in~$U$.
In this representation, the independence
of the functions~$y^1,\ldots,y^k$, $p_2,\ldots,p_k$
translates into the condition that the mapping
\begin{equation*}
(u,[a])\colon \ U\to \bbR^k\times\bbR\bbP^{k-1}
= \bbP\bigl(T^*\bbR^k\bigr)
\end{equation*}
be a submersion.

In fact, the representation~\eqref{eq: P-Dnormalform2}
is more common in treatises on mathematical economics,
where the Pfaff--Darboux theorem plays an important role~\cite{MR1655519}.
Often, the normal forms needed in these models
sometimes come with an additional requirement of \emph{convexity},
i.e., the underlying manifold is $M=\bbR^n$
(or an open domain in~$\bbR^n$),
and one would like to arrange that the functions~$a_i$
be \emph{positive} and the functions~$u^i$ be \emph{strictly convex},
i.e., have positive definite Hessians.\footnote{Sometimes one only requires \emph{weak convexity},
i.e., that the Hessian of $u^i$ be positive definite
on each of its level sets.}
A~useful reference for the role of convexity
in economic models is the book~\cite{ChiapporiEkeland2009}
and the article~\cite{BrowningChiappori1998}.

Now, it turns out that constructing such a \emph{convex} Pfaff--Darboux
representation is not always possible,
which raises the question of how to determine when one exists.
In~\cite{MR2023129}, Ekeland and Nirenberg
were able to provide necessary and sufficient conditions
for a given $1$-form~$\omega\in\Omega^1(\bbR^n)$
to admit a local convex Pfaff--Darboux normal form.
They also constructed examples that showed
that some earlier attempts~\cite{MR1655519, MR1652709}
to find such conditions had been unsuccessful.

In this note, after providing some necessary background,
I prove a generalization of the convex Pfaff--Darboux theorem
of Ekeland and Nirenberg. This treatment has some
notable features that make it of interest for the general problem.

First, the proof of Ekeland and Nirenberg
does not assume the classical Pfaff--Darboux theorem;
instead, it constructs the desired convex representation directly
using the Frobenius theorem, essentially reproving the Pfaff--Darboux theorem
but with the additional convexity condition imposed.
The proof in this article assumes the classical
Pfaff--Darboux theorem, so that the argument can more directly focus
on choosing a Pfaff--Darboux representation
that satisfies the convexity requirements. This results in a shorter
proof, one that also brings the nature of the convexity requirements
more sharply into focus.

Second, the notion of \emph{strict convexity}
turns out to be meaningful on any manifold endowed with
a torsion-free affine connection, and the proof below
covers this more general situation with no extra work.

Third, the proof yields a stronger result,
in that it produces a local convex Pfaff--Darboux representation
of~$\omega$ adapted to any $\omega$-Legendrian foliation
that satisfies a certain geometrically natural positivity condition,
one that is equivalent pointwise to the condition of Ekeland and Nirenberg.

\section{Classical Pfaff--Darboux theorems}
Let $\omega$ be a smooth $1$-form on an $n$-manifold~$M$
that, for some integer $k>0$, satisfies
\begin{equation}\label{eq: omegaPfaff_k}
\text{$\omega\w(\d\omega)^k$ vanishes identically on $M$}
\end{equation}
while
\begin{equation}\label{eq: omega_nondegenerate}
\text{$\omega\w(\d\omega)^{k-1}$ is nowhere vanishing on $M$.}
\end{equation}
The integer~$k{-}1$ is known as the \emph{Pfaff rank}
of~$\omega$~\cite[Chapter II, Section~3]{MR1083148}.
Note that $k\le \tfrac12(n{+}1)$. When $k=\tfrac12(n{+}1)$,
$\omega$ is said to be a \emph{contact form} on~$M$.

When $\omega$ satisfies~\eqref{eq: omegaPfaff_k}
and~\eqref{eq: omega_nondegenerate}, so does $\tilde\omega=f\omega$
for any nonvanishing function~$f$ on~$M$,
since $\tilde\omega\w(\d\tilde\omega)^{r-1} = f^r \omega\w(\d\omega)^{r-1}$
for all integers~$r>0$.

\subsection{Canonical subbundles}
An $\omega$ satisfying~\eqref{eq: omegaPfaff_k}
and~\eqref{eq: omega_nondegenerate} defines a \emph{kernel subbundle}~$K
= \omega^{-1}(0)\subset TM$ of corank~$1$
and a subbundle $A\subset K$ of corank $2(k{-}1)$ in~$K$
by the rule that, for each $m\in M$,
\begin{equation*}
A_m = \{ v\in K_m\mid \mathrm{d}\omega(v,w) = 0,\, \forall w\in K_m \}.
\end{equation*}
Replacing~$\omega$ by $\tilde\omega= f \omega$ for any nonvanishing
function~$f$ does not change~$K$ or $A$.

If $k>1$, then $K\subset TM$ is not an integrable plane field,
but the subbundle $A\subset TM$ is always integrable,
since it is the Cauchy characteristic plane field
of the differential ideal~$\cI$ generated by~$\omega$
(see~\cite[Chapter II, Proposition~2.1]{MR1083148}).
In the contact case, i.e., when $n=2k{-}1$ (which is,
in some sense, generic), one has~$A=(0)$.

There is a nondegenerate, skew-symmetric bilinear pairing
$B_\omega \colon K{/}A\times K{/}A\to\bbR$
defined by
\[
B_\omega\bigl(v{+}A_m,w{+}A_m\bigr) = \mathrm{d}\omega(v,w),
\]
when $v,w\in K_m$. It satisfies $B_{f\omega} = fB_\omega$
for any nonvanishing function~$f$ on~$M$.

Note that any subspace~$W\subset T_m$ on which both $\omega$ and $\d\omega$
vanish, must, first of all, satisfy $W\subset K_m$
(since $\omega$ vanishes on~$W$), and, second, must have codimension
at least $k{-}1$ in $K_m$, since $\d\omega$, as a skew-symmetric form
on $K_m$, has Pfaff rank~$k{-}1$.
Moreover, if $W$ does have codimension $k{-}1$ in $K_m$,
then it must contain $A_m$, so that $W/A_m$ is a null subspace of $B_\omega$.

\subsection{Legendrian submanifolds and Grassmannians}
Any submanifold~$L\subset M$ that satisfies~$L^*\omega=0$,
i.e., an \emph{integral manifold of~$\omega$},
must also satisfy $L^*\d\omega=0$
and hence, by the above linear algebra discussion,
must have codimension at least~$k$ in~$M$.
When $L\subset M$ is an integral manifold of~$\omega$ of codimension~$k$,
it is said to be an \emph{$\omega$-Legendrian submanifold}.

In particular, if $L\subset M$ is $\omega$-Legendrian,
then, for each $m\in L$, the tangent space $T_mL$
satisfies $A_m\subset T_mL\subset K_m$
while $B_\omega$ vanishes identically on $T_mL/A_m\subset K_m/A_m$.

This motivates defining
the \emph{Legendrian Grassmannian}~$\Leg_m(\omega)\subset\Gr^k(T_mM)$
to be the set of subspaces~$W\subset K_m$
that have codimension~$k$ in $T_mM$
and on which both~$\omega$ and~$\d\omega$ vanish.
By the above remarks, it follows that $\Leg_m(\omega)$
can be canonically identified with the
\emph{Lagrangian Grassmannian}~$\Lag(K_m/A_m)\subset\Gr^{k-1}(K_m/A_m)$
consisting of the $(k{-}1)$-dimensional subspaces of~$K_m/A_m$
on which $B_\omega$ vanishes.
Hence, $\Leg_m(\omega)$ is naturally a smooth manifold
of dimension~$\tfrac12k(k{-}1)$.
Moreover, the (disjoint) union
\[
\Leg(\omega) = \bigcup_{m\in M} \Leg_m(\omega)\subset \Gr^k(TM)
\]
is a smooth subbundle,
and $\Leg(f\omega)=\Leg(\omega)$ for all nonvanishing~$f$.

\subsection{A local normal form}
One version of the \emph{Pfaff--Darboux theorem}
\cite[Chapter II, Theorem~3.1]{MR1083148} states that,
when ${\omega\in\Omega^1(M)}$ satisfies~\eqref{eq: omegaPfaff_k}
and~\eqref{eq: omega_nondegenerate},
each $m\in M$
has an open neighborhood~$U\subset M$
on which there exist smooth functions~$u^1,\ldots,u^k$
and $a_1,\ldots,a_k$ (with not all $a_i$ simultaneously vanishing)
so that
\begin{equation}\label{eq: PfaffDarboux_normalform}
U^*\omega = a_1\,\d u^1 + \cdots + a_k\,\d u^k .
\end{equation}
Moreover,
the mapping~$(u,[a])\colon U\to\bbR^k\times\mathbb{RP}^{k-1}$
is a submersion. (In fact, the kernel subbundle of the differential
of this mapping is the restriction of~$A$ to~$U$.)

Conversely, the existence of functions $u^i$ and $a_i$ for~$1\le i\le k$
on an open set~$U\subset M$ satisfying~\eqref{eq: PfaffDarboux_normalform}
with the $a_i$ not all simultaneously vanishing and having the
property that $(u,[a]) \colon U\to\bbR^k\times\mathbb{RP}^{k-1}$
be a submersion implies that both~\eqref{eq: omegaPfaff_k}
and~\eqref{eq: omega_nondegenerate} hold on~$U$.

\subsection{Geometry of the normal form}
It will be useful to have a geometric interpretation
of the Pfaff--Darboux theorem.
Now, in the representation~\eqref{eq: PfaffDarboux_normalform},
the functions~$u^i$ have independent differentials,
i.e., $\d u^1\w\cdots\w\d u^k$ does not vanish on~$U$.
Consequently, the simultaneous level sets of the functions~$u^i$
define a foliation~$\cL$ of~$U\subset M$
by $\omega$-Legendrian submanifolds, i.e., an $\omega$-Legendrian foliation.

Conversely, given an $\omega$-Legendrian foliation~$\cL$
on an open subset~$V\subset M$, each point~$m\in V$ will
have an open neighborhood~$U\subset V$ in which the leaves of~$\cL$
are the fibers of a submersion~${u = (u^i) \colon U\to\bbR^k}$.
Since~$\omega$ vanishes when pulled back to any fiber of~$u$,
it follows that there exists a mapping~$a = (a_i) \colon U\to\bbR^k$
such that $U^*\omega = a_1\,\d u^1 + \cdots + a_k\,\d u^k$.

Thus, a geometric interpretation of the Pfaff--Darboux theorem
is the statement that,
when $\omega\in\Omega^1(M)$ satisfies~\eqref{eq: omegaPfaff_k}
and~\eqref{eq: omega_nondegenerate}, each point~$m\in M$
has an open neighborhood~$U\subset M$
on which there exists an $\omega$-Legendrian foliation.

\subsection{Variants and extensions}
There are a number of variants and extensions
of the classical Pfaff--Darboux theorem
that can all be seen to be equivalent to the above versions
by elementary arguments~\cite[Chapter II, Section~3]{MR1083148}.
In this article, two such variants will be important.
For convenience of reference, they will be stated as propositions.

\begin{Proposition}\label{prop: Legfromtangent}
Suppose that $\omega\in\Omega^1(M)$ satisfies~\eqref{eq: omegaPfaff_k}
and~\eqref{eq: omega_nondegenerate}. Then for each $m\in M$
and $W\in\Leg_m(\omega)$,
there exists a $\omega$-Legendrian submanifold~$L\subset M$
such that $m\in L$ and $W = T_mL$.
\end{Proposition}

\begin{Proposition}\label{prop: folfromleg}
Suppose that $\omega\in\Omega^1(M)$ satisfies~\eqref{eq: omegaPfaff_k}
and~\eqref{eq: omega_nondegenerate} and that $L\subset M$
is an embedded $\omega$-Legendrian submanifold. Then each $m\in L$
has an open neighborhood~$U\subset M$ on which there exists
an $\omega$-Legendrian foliation~$\cL$ with the property that $L\cap U$
is a leaf of $\cL$.
\end{Proposition}

\section{Convexity and affine manifolds}

\subsection{Classical convexity}
When~$M=\bbR^n$, there is a notion of \emph{strict convexity}
of a function~$u$, which is the condition
that the Hessian quadratic form~$H(u)$ be positive definite, where
\begin{equation}\label{eq: flat_Hessian}
H(u) = {\frac{\partial^2 u}{\partial x^i\partial x^j}}
\,\d x^i\otimes\d x^j
\end{equation}
and where $x^1,\ldots,x^n$
are the usual affine linear coordinates in~$\bbR^n$.
Note that strict convexity is an affine-invariant notion on~$\bbR^n$.

Motivated by applications in economics,
Ekeland and Nirenberg~\cite{MR2023129}
asked what further conditions one must impose on an~$\omega\in\Omega^1(\bbR^n)$
satisfying~\eqref{eq: omegaPfaff_k} and~\eqref{eq: omega_nondegenerate}
in order to know that one can choose the functions~$u^j$
and $a_j$ in the representation~\eqref{eq: PfaffDarboux_normalform}
so that the $u^j$ be strictly convex and the $a_j$ be positive.
It is not hard to show, by example, that \emph{some}
further condition on~$\omega$ is necessary
to guarantee the existence of such a convex representation.
(See the discussion at the beginning of Section~\ref{ssec: poscond}.)

They showed that two earlier articles~\cite{MR1655519,MR1652709}
claiming to provide such necessary and sufficient conditions were flawed
(indeed, they exhibited counterexamples to the claims of these articles)
and then produced their own condition,
which they showed to be necessary and sufficient.

In this note, I will show that their main result,
properly formulated, holds good on an $n$-manifold~$M$
endowed with a torsion-free affine connection,
not just on~$\bbR^n$ endowed with the (flat) affine connection
it inherits as a vector space.

\subsection{Affine connections and convexity}

Let~$\nabla$ be a torsion-free affine connection on an $n$-manifold~$M^n$,
i.e., $\nabla$ is a first-order, linear differential operator
\begin{equation*}
\nabla \colon \ \Omega^1(M)\to \Omega^1(M)\otimes \Omega^1(M)
\end{equation*}
that obeys the Leibnitz rule
\begin{equation*}
\nabla (f\eta) = \d f\otimes\eta + f \nabla(\eta)
\end{equation*}
for all smooth functions~$f$ on~$M$ and smooth $1$-forms~$\eta$ on~$M$.
The assumption that~$\nabla$ be torsion-free
is the condition that the associated (second-order)
\emph{Hessian operator}~$H(u)=\nabla(\d u)$
be a symmetric $(0,2)$-tensor for each smooth function~$u$ on~$M$.

A smooth function~$u$ on~$M$ is said to be {\it strictly $\nabla$-convex}
if, as a quadratic form,~$H(u)$ is positive definite at every point of~$M$.

When~$M=\bbR^n$ and~$\nabla$ is the standard (flat) connection,
satisfying~$\nabla(\d x^i) = 0$ for all of the coordinate functions~$x^i$,
then~$H(u)$ is the usual Hessian tensor~\eqref{eq: flat_Hessian},
and this notion of convexity is simply the classical one.

In the more general case, when $x = (x^i)\colon U\to \bbR^n$ is a local
coordinate chart, one has
\begin{equation*}
H\big(x^k\big) = \nabla\big(\d x^k\big) = \Gamma^k_{ij}\,\d x^i\otimes\d x^j,
\end{equation*}
where $\Gamma^k_{ij}=\Gamma^k_{ji}\in C^\infty(U)$
are the coefficients of the connection $\nabla$ relative
to the coordinate chart~$x=(x^i)$.
The general coordinate formula for~$H$ then becomes
\begin{equation*}
H(u) = \left(\frac{\partial^2 u}{\partial x^i\partial x^j}
+ \Gamma^k_{ij} \frac{\partial u}{\partial x^k}\right)
\,\d x^i\otimes\d x^j.
\end{equation*}
Thus, $\nabla$-convexity of $u$ is expressible in terms of a condition
on the $2$-jet of~$u$, slightly more general than the condition for
classical convexity.

Adopting the usual conventions
\begin{equation*}
\alpha\w\beta
= \tfrac12(\alpha\otimes\beta-\beta\otimes\alpha), \qquad
\alpha\circ\beta
= \tfrac12(\alpha\otimes\beta+\beta\otimes\alpha),
\end{equation*}
one sees that, for a $1$-form~$\omega$ of the form
\begin{equation}\label{eq: P-Drepresentation}
\omega = a_1\,\d u^1 + \cdots + a_k\,\d u^k,
\end{equation}
one has (using the summation convention)
\begin{align*}
\nabla\omega
&= \d a_i\otimes\d u^i + a_i \nabla(u^i) = \d a_i\otimes\d u^i + a_i H(u^i)\\
&= \d a_i\w\d u^i + \d a_i\circ\d u^i + a_i H(u^i)\\
&= \d\omega + \Ss^\nabla\omega,
\end{align*}
where I have introduced
the notation~$\Ss^\nabla\omega$ to denote the {\it symmetrization\/}
of~$\nabla(\omega)$. Thus,~$\Ss^\nabla\omega
= \nabla\omega- \d\omega$
is a well-defined quadratic form on~$M$.
(Of course, the linear, first-order differential operator~$\Ss^\nabla$
depends on~$\nabla$.)

\subsection{A positivity condition}\label{ssec: poscond}
The quadratic form $S^\nabla\omega$ provides some insight
into the question of whether a $\nabla$-convex Pfaff--Darboux
representation of $\omega$ is possible.

\begin{Proposition}\label{prop: ness cond}
Suppose that $\omega\in\Omega^1(M)$ satisfies~\eqref{eq: omegaPfaff_k}
and~\eqref{eq: omega_nondegenerate}. If there exist positive
functions $a_i$ and $\nabla$-convex functions $u^i$ for $1\le i\le k$
such that~\eqref{eq: P-Drepresentation} holds, then $S^\nabla\omega$
is positive definite on the leaves of the foliation $\cL$ defined by
$\d u^1 = \d u^2 = \cdots = \d u^k = 0$.
\end{Proposition}

\begin{proof}
Since $\Ss^\nabla\omega =\d a_i\circ \d u^i + a_i H(u^i)$,
it follows that, when restricted to the plane field ${N\subset TM}$
defined by~$\d u^1 = \d u^2 = \cdots = \d u^k = 0$,
the terms~$\d a_i\circ \d u^i$ in~$\Ss^\nabla\omega$
vanish. Thus, $S^\nabla\omega = a_i H(u^i)$
as quadratic forms on $N$. By the positivity of the $a_i$
and the $\nabla$-convexity of the $u^i$, it follows that
$S^\nabla\omega$ is positive definite on~$N$.\footnote{It is worth pointing out that the same conclusion
about the positive definiteness of $\Ss^\nabla\omega$ on the leaves of~$\cL$
would have followed if one had merely assumed that each $u^i$
be only `strictly $\nabla$-\emph{quasi}-convex',
i.e., that $\d u^i$ be nonvanishing and $H(u^i)$ be positive definite
when restricted to the hyperplane field~$\d u^i = 0$.
Compare~\cite[Lemma~1]{MR2023129},
and the preceding discussion about their Problem~2.}
\end{proof}

This proposition provides necessary condition for the existence
of a convex Pfaff--Darboux representation.

\begin{Example}[an obstructed example]\label{ex: Obstruction}
Let $M=\mathbb{R}^n$ with standard coordinates $x=(x^i)$,
and let $\nabla$ be the (flat, torsion-free) connection
such that $\nabla(\d x^i)=0$ for $1\le i\le n$.
Let $c_i$, $f_{ij}=-f_{ji}$, and $g_{ij}=g_{ji}$
be constants and consider the $1$-form
\[
\omega = \bigl(c_i + (f_{ij}+g_{ij})x^j)\,\d x^i.
\]
Then $\d\omega= f_{ij}\,\d x^j\w\d x^i = - f_{ij}\,\d x^i\w\d x^j$
and $\Ss^\nabla\omega = g_{ij}\,\d x^i\circ\d x^j$.

Now assume that the skew-symmetric matrix $f = (f_{ij})$
has rank $2(k{-}1)<n$, so that $(\d\omega)^{k-1}\not=0$
but $(\d\omega)^{k}=0$.
Then, for generic choice of the constants $c_i$,
$\omega\w(\d\omega)^{k-1}$ will be nonvanishing on an open neighborhood $U\subset\mathbb{R}^n$ of the origin $0\in\mathbb{R}^n$,
in which case, $\omega$ will satisfy the hypotheses~\eqref{eq: omegaPfaff_k}
and~\eqref{eq: omega_nondegenerate} on~$U$.

If the symmetric matrix $g = (g_{ij})$
does not have at least $n{-}k$ positive eigenvalues,
then~$\Ss^\nabla\omega$ cannot be positive definite on
any $(n{-}k)$-dimensional subbundle $N\subset TU$, and, hence,
by Proposition~\ref{prop: ness cond}, $\omega$ cannot have a convex Pfaff--Darboux
representation in any open subset of $U$.
\end{Example}

It turns out that this necessary condition
for a local `convex' Pfaff--Darboux representation compatible
with an $\omega$-Legendrian foliation~$\cL$ is also sufficient.

\begin{Theorem}\label{thm: main}
Suppose $\nabla$ be a torsion-free affine
connection on~$M$, that $\omega\in\Omega^1(M)$
satisfy~\eqref{eq: omegaPfaff_k} and~\eqref{eq: omega_nondegenerate}
for some~$k>0$, and that $\cL$ be an $\omega$-Legendrian
foliation on~$M$ with the property that $\Ss^\nabla\omega$
pulls back to each leaf of~$\cL$ to be positive definite.
Then each $m\in M$ has an open neighborhood~$U\subset M$
on which there exist strictly $\nabla$-convex functions~$u^1,\ldots,u^k$
that are constant on the leaves of~$\cL$ in~$U$
and positive functions~$a_1,\ldots, a_k$ such that
\begin{equation*}
U^*\omega = a_1\,\d u^1 + \cdots + a_k\,\d u^k.
\end{equation*}
\end{Theorem}

Before giving the proof of Theorem~\ref{thm: main},
I will state one of its corollaries, so that it can be compared
with the main result of Ekeland and Nirenberg~\cite[Theorem 1]{MR2023129}.

First, some useful terminology.
As always, assume that $\omega$ satisfies
\eqref{eq: omegaPfaff_k} and~\eqref{eq: omega_nondegenerate} for some~${k>0}$.

\begin{Definition} An $\omega$-Legendrian subspace~$W\subset T_mM$
is \emph{$\nabla$-positive for~$\omega$}
if the restriction of the quadratic form~$\Ss^\nabla\omega$ to~$W$
is positive definite.
\end{Definition}

Let $\Leg^+(\omega,\nabla)\subseteq\Leg(\omega)$
denote the set of $\omega$-Legendrian subspaces
that are $\nabla$-positive for~$\omega$.
Then $\Leg^+(\omega,\nabla)$ is a (possibly empty) open subset
of~$\Leg(\omega)$. Consequently, the set of points~$m\in M$
for which there exists a $\nabla$-positive, $\omega$-Legendrian
subspace~$W\subset T_mM$ is an open subset of~$M$.
Also, note that, since such a $W$ contains~$A_m$,
it follows that $\Ss^\nabla\omega$ must be positive definite on~$A_m$.

\begin{Corollary}\label{cor: main}
Suppose that $\nabla$ be a torsion-free affine
connection on~$M$, that $\omega\in\Omega^1(M)$
satisfy~\eqref{eq: omegaPfaff_k} and~\eqref{eq: omega_nondegenerate}
for some~$k>0$, and that there exist a $W\in\Leg^+(\omega,\nabla)$
with ${W\subset T_mM}$. Then $m\in M$ has an open neighborhood ${U\subset M}$
on which there exist strictly $\nabla$-convex functions~$u^1,\ldots,u^k$
and positive functions~$a_1,\ldots,a_k$ such that
\begin{equation*}
U^*\omega = a_1\,\d u^1 + \cdots + a_k\,\d u^k.
\end{equation*}
\end{Corollary}

The proof of Corollary~\ref{cor: main} follows
by applying Propositions~\ref{prop: Legfromtangent}
and~\ref{prop: folfromleg} to produce an $\omega$-Legendrian foliation~$\cL$
on an open neighborhood~$V$ of~$m$ whose leaf through~$m$ has~$W$
as a~tangent space. Since $\Ss^\nabla\omega$ is positive definite on~$W$,
it follows that it is positive definite on all the tangent spaces
to the leaves of~$\cL$
in some (possibly) smaller $m$-neighborhood~$V'\subset V$.
Now apply Theorem~\ref{thm: main} to~$\cL$ on~$V'$.

\begin{Remark}
In the special case in which~$M=\bbR^n$ and $\nabla$ is the flat
connection satisfying~$\nabla(\d x^i)=0$ for $x^i$ the standard
coordinates on~$\bbR^n$, Corollary~\ref{cor: main}
simply becomes Theorem~1 of Ekeland and Nirenberg~\cite{MR2023129},
since their Condition~3 turns out to be equivalent
to the existence of a~$W\in\Leg^+(\omega,\nabla)$ with $W\subset T_mM$
in this case.
\end{Remark}

\begin{proof}[Proof of Theorem~\ref{thm: main}.]
There exists an $m$-neighborhood $V_0\subset M$
on which there exist smooth functions~$y^1,\ldots,y^k$
vanishing at~$m$
so that the leaves of $\d y^1 =\cdots =\d y^k = 0$
are intersections of the leaves of~$\cL$ with~$V_0$
as well as functions~$p_2,\ldots,p_k$, also vanishing
at~$m$, and a nonvanishing function~$a$ so that
\[
V_0^*\omega
= a\bigl(\d y^1 + p_2\,\d y^2 +\cdots + p_k\,\d y^k\bigr)
\]
By reversing the signs of~$a$ and the~$y^i$,
if necessary, one can assume that~$a(m)>0$.
Let~$W\subset T_mM$ be the tangent to the leaf of~$\cL$
that passes through~$m$, so that $W$ is the common kernel
of the~$\d y^i$ evaluated at~$m$.

Set~$\bar\omega = a^{-1}\omega$ and note that,
since~$\d\bar\omega \equiv a^{-1}\,\d\omega \bmod \omega$,
it follows that $\cL$ is also $\bar\omega$-Legendrian.
Moreover, since
\[
\Ss^\nabla\bar\omega = \d\big(a^{-1}\big)\circ\omega + a^{-1} \Ss^\nabla\omega,
\]
it follows that the tangent spaces of~$\cL$
(which, of course, satisfy $\omega=0$)
are also $\nabla$-positive for~$\bar\omega$.
Since~$\omega = a \bar\omega$ and $a>0$,
finding the desired convex representation for $\bar\omega$
will also yield one for~$\omega$.
Thus, it suffices to prove the theorem with~$\bar\omega$
in the place of~$\omega$, i.e., to assume that~$a=1$,
so I will do that from now on. Thus,
\[
V_0^*\omega = \d y^1 + p_2\,\d y^2 +\cdots + p_k\,\d y^k.
\]
Since~$\omega\w(\d\omega)^{k-1}\not=0$, the functions
$y^1,\ldots,y^k$, $p_2,\ldots,p_k$ have linearly independent differentials
at~$m$.

Restricting to~$V_0$, i.e., setting $M=V_0$, one has
\[
\Ss^\nabla\omega = H\big(y^1\big) + \d p_2\circ \d y^2 + \cdots + \d p_k\circ \d y^k
 + p_2 H\big(y^2\big) + \cdots + p_k H\big(y^k\big).
\]
Since the~$p_j$ vanish at~$m$,
it follows that, when restricted to~$W\subset T_mM$,
the two quadratic forms~$H\big(y^1\big)$ and $\Ss^\nabla\omega$ are equal.
Thus~$H\big(y^1\big)$ is positive definite on~$W$,
and so there is a~constant~$c>0$
so that~\smash{$H\big(y^1\big) + c \big(\d y^2\big)^2 + \cdots + c \big(\d y^k\big)^2$}
is positive definite on~$K_m = \{v\in T_mM\ \vrule\ \d y^1(v) = 0\}$.
Writing
\[
\omega = \d\bigl(y^1 + \tfrac12c \big(y^2\big)^2 + \cdots + \tfrac12c \big(y^k\big)^2\bigr)
 + \big(p_2-c y^2\big)\,\d y^2 + \cdots + \big(p_k-c y^k\big)\,\d y^k
\]
and observing that
\begin{gather*}
H\big(y^1 + \tfrac12c \big(y^2\big)^2 + \cdots + \tfrac12c \big(y^k\big)^2\big)\\
\qquad = H\big(y^1\big) + c \big(\d y^2\big)^2 + \cdots + c \big(\d y^k\big)^2
 + c y^2 H\big(y^2\big) + \cdots + c y^k H\big(y^k\big)
\end{gather*}
shows that, setting
\[
\bar y^1 = y^1 + \tfrac12c \big(y^2\big)^2 + \cdots + \tfrac12c \big(y^k\big)^2,
\qquad
\bar y^i = y^i,
\qquad\text{and}\qquad
\bar p_i = p_i - c y^i,
\]
one has
$\omega = \d\bar y^1
 + \bar p_2\,\d \bar y^2 + \cdots + \bar p_k\,\d \bar y^k$.

Thus, one could have chosen the functions~$y^1,\ldots,y^k$, $p_2,\ldots,p_k$
with~$H\big(y^1\big)$ being positive definite on the hyperplane~$K_m$.
Assume now that this was done.

It still needs to be shown that one can choose
the functions~$y^1,\ldots,y^k$, $p_2,\ldots,p_k$
with $H\big(y^1\big)$ being positive definite on all of $T_mM$,
not just on $K_m$, which is the kernel of $\d y^1$ at~$m$.
To do this, note that, if~$\phi$ is any smooth function on a neighborhood
of the origin in~$\bbR$, then
\[
H\bigl(\phi\big(y^1\big)\bigr) = \phi'\big(y^1\big) H\big(y^1\big) + \phi''\big(y^1\big) \big(\d y^1\big)^2.
\]
Hence, by choosing a~$\phi$ with $\phi(0)=0$, $\phi'(0)=1$
and $\phi''(0)>0$ sufficiently large,
I can arrange that~$\phi\big(y^1\big)$ be strictly $\nabla$-convex at~$m$.
Since
\[
\omega = \frac{1}{\phi'(y_1)} \big(\d\bigl(\phi\big(y^1\big)\bigr)
    + \phi'(y_1)p_2\,\d y^2 +\cdots + \phi'(y_1)p_k\,\d y^k\big),
\]
one sees that the functions
$\bar y^1,\ldots,\bar y^k$, $\bar p_2,\ldots,\bar p_k$, where
\[
\bar y^1 = \phi\big(y^1\big),
\qquad\text{and}\qquad
\bar y^i = y^i,
\qquad
\bar p_i = \phi'(y_1)p_i,\qquad 2\le i\le k,
\]
(with $a = 1/\phi'(y_1) >0$), give a Pfaff--Darboux representation
for~$\omega$ that is compatible with the foliation~$\cL$ and
for which $\bar y^1$ is strictly $\nabla$-convex.\footnote{This is the same idea that Ekeland and Nirenberg \cite{MR2023129}
used in their generalization of their Theorem~1
to cover the quasi-convex case.}

Thus, one can assume henceforth that,
on an open $m$-neighborhood~$V_1\subset M$,
one has a~representation of the form
\[
\omega = \d y^1 + p_2\,\d y^2 +\cdots + p_k\,\d y^k,
\]
where the functions~$y^1,\ldots,y^k,p_2,\ldots,p_k\in C^\infty(V_1)$
all vanish at~$m$,
the equations $\d y^i=0$ define the tangents to the leaves of~$\cL$ in~$V_1$,
and~$y^1$ is strictly~$\nabla$-convex.

Under these assumptions, there is a constant~$b>0$ sufficiently large
so that $H\big(y^i+by^1\big) = H\big(y^i\big)+b H\big(y^1\big)$ is positive definite at~$m$
for~$2\le i\le k$.
Thus, writing
\[
\omega = \bigl(1-b(p_2{+}\cdots{+}p_k)\bigr)\,\d y^1 + p_2\,\d\big(y^2+b y^1\big)
   + \cdots + p_k\,\d\big(y^k+b y^1\big),
\]
it follows that I can,
after restricting to an $m$-neighborhood $V_2\subset V_1$ on which
the function $a = \bigl(1-b(p_2{+}\cdots{+}p_k)\bigr)$ is positive,
dividing by~$a>0$,
and replacing~$y^j$ by $y^j+by^1$ and~$p_j$ by~$p_j/a$ for $2\le j\le k$,
assume that I have a representation
\[
\omega = \d y^1 + p_2\,\d y^2 +\cdots + p_k\,\d y^k,
\]
in which all of the~$H\big(y^j\big)$ are positive definite at~$m$,
i.e., the $y^j$ are strictly $\nabla$-convex on some neighborhood of~$m$
and the $p_i$ all vanish at~$m$.

Finally, for $\varepsilon>0$ and sufficiently small, write
\[
\omega = \d\bigl(y^1-\varepsilon\big(y^2{+}\cdots{+}y^k\big)\bigr)
   + (p_2 + \varepsilon)\,\d y^2
   + \cdots + (p_k + \varepsilon)\,\d y^k.
\]
Then, setting~$u^1 = y^1-\varepsilon\big(y^2+\cdots+y^k\big)$ and~$u^j = y^j$
for~$j>1$ and setting~$a_1 = 1$ and~${a_j = \varepsilon + p_j}$ for~$j>1$,
one achieves the desired convex Pfaff--Darboux representation
on an open $m$-neigh\-bor\-hood~$U\subset V_2$.
\end{proof}

\begin{Example}[dependence on $\nabla$]
Let $\omega$ be defined on $\mathbb{R}^n$ as in Example~\ref{ex: Obstruction}
with the constants~$c_i$,~$f_{ij}$ and $g_{ij}$ as specified there
so that $\omega\w (\d\omega)^k=0$ and $\omega\w(\d\omega)^{k-1}$
is non-vanishing on an open neighborhood $U\subset\mathbb{R}^n$ of $0\in\mathbb{R}^n$. Now, however, suppose that $\nabla$ is the
torsion-free connection that satisfies $\nabla(\d x^i) = \Gamma^i_{jk}\,\d x^j\otimes \d x^k$, where $\Gamma^i_{jk}=\Gamma^i_{kj}$ are constants.
Then calculation yields
\[
\Ss^\nabla\omega
= \bigl(g_{ij}+\bigl(c_k+(f_{kl}+g_{kl})x^l\bigr)\Gamma^k_{ij}\bigr)
\,\d x^i\circ\d x^j.
\]
Consequently, at the origin $x=0$, one finds $\big(\Ss^\nabla\omega\big)_0
= \bigl(g_{ij}+c_k\Gamma^k_{ij}\bigr)\,\d x^i\circ\d x^j$. It follows
that, no matter the value of $g = (g_{ij})$, one can always choose
$\Gamma^i_{jk}$ so that $\big(\Ss^\nabla\omega\big)_0$
be positive definite.
Thus, by Theorem~\ref{thm: main} there always exists a
torsion-free connection~$\nabla$
for which a $\nabla$-convex Pfaff--Darboux representation of $\omega$
exists on a neighborhood of $0\in\mathbb{R}^n$.

This highlights the significance of the choice of background connection
for the convex representability of a given~$\omega$. In turn,
this makes clear how the choice of coordinates in which a given
problem is described affects the existence of convex representability.
\end{Example}

\begin{Remark}[global considerations]
While Theorem~\ref{thm: main} gives necessary and sufficient conditions
for the existence of local $\nabla$-convex Pfaff--Darboux representations,
for applications one would like to know something about how large an open
set in the model~$M$ one can cover with such a representation,
and this seems to be a subtle problem.

Even in the simplest case of a $3$-manifold~$M$
endowed with a contact $1$-form~$\omega$
and a torsion-free affine connection~$\nabla$
for which $\Ss^\nabla\omega$ is positive definite
on the $2$-plane bundle~$K\subset TM$,
it is not clear how to characterize the domains~$U\subset M$
that support a $\nabla$-convex Pfaff--Darboux representation for~$\omega$.
\end{Remark}

\begin{Example}[global non-existence]
Let $S^3\subset\mathbb{R}^4\simeq\mathbb{H}$
be the unit sphere regarded as the Lie group of unit quaternions.
Let $\omega_i$ for~$i=1,2,3$ be a basis for the left-invariant $1$-forms
on $S^3$, which obey the formulae $\d\omega_i = \epsilon_{ijk} \omega_j\w\omega_k$ where $\epsilon$ is the fully skew-symmetric symbol
satisfying $\epsilon_{123}=1$. In particular, $\omega_i\w(\d\omega_i)
= 2 \omega_1\w\omega_2\w\omega_3\not=0$ but $\omega_i\w(\d\omega_i)^2 = 0$,
so that $\omega_i$ satisfies the usual hypotheses with $k=1$, i.e.,
it is a contact form on $S^3$.

Now let $\Gamma_{ijk}=\Gamma_{ikj}$ be constants for $1\le i,j,k\le 3$
and let $\nabla$ be the affine connection (necessarily torsion-free) that satisfies
\[
\nabla(\omega_i) = (\epsilon_{ijk}+\Gamma_{ijk}) \omega_j\otimes\omega_k .
\]
Then
\[
S^\nabla(\omega_i) = \Gamma_{ijk} \omega_j\circ\omega_k.
\]

Note that, if we choose $\Gamma_{ijk}=0$, so that $\nabla$
is the Levi-Civita connection of the constant curvature metric
$g = {\omega_1}^2+{\omega_2}^2+{\omega_3}^2$ on $S^3$,
then $S^\nabla(\omega_i)\equiv0$, so $\omega_i$ cannot
have a $\nabla$-convex Pfaff--Darboux representation
on any open set in $S^3$.

However, with appropriate choice of constants $\Gamma_{ijk}$,
we could arrange that
$S^\nabla(\omega_i)$ be positive definite for each $i$. In this case,
by~\cite[Theorem~1]{MR2023129}, each point $m\in S^3$ will have an open neighborhood
on which $\omega_i$ has a $\nabla$-convex Pfaff--Darboux representation.

However, because $S^3$ is compact, any smooth function $u$ on $S^3$
must have a local maximum, and $H(u)$ will be negative semidefinite
there, independent of the choice of $\nabla$. Thus, there cannot be
any global functions $u$ on $S^3$ that are $\nabla$-convex.
\emph{A fortiori,}
there cannot be a global $\nabla$-convex Pfaff--Darboux representation
for $\omega_i$ for any $i$.
\end{Example}

\subsection*{Acknowledgements}

This article is dedicated to Jean-Pierre Bourguignon,
with much admiration, on the occasion of his 75th birthday.
Thanks to Duke University for its support
via a research grant and to the National Science Foundation
for its support via DMS-9870164
(during which most of the research for this article was done)
and DMS-1359583 (during which this article was written).

\pdfbookmark[1]{References}{ref}
\LastPageEnding

\end{document}